\documentclass[a4paper, reqno]{amsart}

\usepackage[all]{xy}
\usepackage[english]{babel}
\usepackage[utf8]{inputenc}
\usepackage{amssymb,amsmath,amsthm}
\usepackage{amsfonts}
\usepackage{graphicx}
\usepackage{psfrag}
\usepackage[dvipsnames]{xcolor}
\usepackage[hidelinks]{hyperref}
\usepackage{csquotes}
\usepackage[alphabetic, msc-links, initials]{amsrefs}
\usepackage{tikz-cd}
\usepackage{esint}

\makeatletter
\@namedef{subjclassname@2020}{\textup{2020} Mathematics Subject Classification}
\makeatother

\theoremstyle{plain}
\newtheorem{theorem}{Theorem}[section]

\theoremstyle{definition}
\newtheorem{definition}[theorem]{Definition}

\newtheorem*{question}{Question}
\newtheorem{example}[theorem]{Example}

\theoremstyle{remark}
\newtheorem{remark}[theorem]{Remark}

\theoremstyle{plain}

\newcommand{\Hyp}{\mathbb{H}}

\newcommand{\R}{\mathbb{R}}

\newcommand{\C}{\mathbb{C}}

\newcommand{\N}{\mathbb{N}}

\newcommand{\scal}{\operatorname{scal}}
\newcommand{\ric}{\operatorname{Ric}}

\newcommand{\volume}{\operatorname{vol}}

\newcommand{\dv}{\text{ }dV}

\newcounter{mnotecount}[section]

\allowdisplaybreaks

\title[Ricci flow and the scalar curvature rigidity of Einstein manifolds]{Ricci flow and the scalar curvature rigidity of Einstein manifolds}

\author{Klaus Kröncke}
\address{Institutionen för Matematik, Kungliga Tekniska Högskolan, 100 44 Stockholm, Sweden}
\email{kroncke@kth.se}


\begin{document}
\hbadness=100000
\vbadness=100000

\begin{abstract}
We review recent results relating linear stability to dynamical stability and the scalar curvature rigidity of Einstein manifolds. We discuss closed and open Einstein manifolds as well as complete noncompact Einstein manifolds which are asymptotically locally Euclidean and asymptotically hyperbolic. For these classes, the relation to the positive mass theorem will also be explained.
\end{abstract}

\maketitle

\begin{sloppypar}

\section{Introduction}
A Riemannian manifold $(M,g)$ is called \textit{Einstein}, if $\ric_g=\mu g$ for some $\mu\in\R$ which we call the \textit{Einstein constant}. For Einstein manifolds, there are some natural notions of stability and rigidity which are closely related. For example, one may ask if $g$ is stable in the sense of being a local extremum of some natural geometric functionals on the space of metrics or \emph{dynamically stable} as a stationary point of an appropriately normalized variant of the Ricci-flow. Such a variant could for example be the volume-normalized Ricci flow
\begin{equation*}
\partial_tg_t=-2\ric_g+\left(\frac{2}{n\cdot \mathrm{vol}(M,g)}\int_M\scal_{g}\dv_{g}\right) g,
\end{equation*}
if the manifold is closed. If it is of infinite volume, one may instead work with the flow
\begin{equation*}
\partial_tg=-2\ric_{g}+2\mu g.
\end{equation*}
One may also ask if an Einstein metric $g$ is \textit{rigid} in the sense that it is an isolated point in the moduli space of Einstein metrics, i.e.\ all other nearby Einstein metrics are up to scaling isometric to $g$. Another rigidity question is \textit{scalar curvature rigidity}, which asks whether $g$ can be deformed on a compact subset to a metric of larger scalar curvature.
These questions are closely related to the spectrum of the Einstein operator $\Delta_E$ acting on symmetric $2$-tensors, which appears in the linearization of the Einstein condition.
It is given by
\begin{equation*}
\Delta_Eh=\nabla^*\nabla h-2\mathring{R}h,\qquad (\mathring{R}h)_{ij}:=R_{iklj}h^{kl}.
\end{equation*}
It is very hard to obtain spectral estimates for this operator. Some results are summarized in an earlier overview article by the author \cite{Kro16}. A bit more can be said if the Einstein manifold has special holonomy, or is a symmetric or a homogenous space. In fact, $\Delta_E$ has recently been intensively investigated on these spaces and these results are discussed in another overview article by Schwahn and Semmelmann in the same volume \cite{SS25}. 
Here, we will mainly focus on discussing the progress which recently has been made in understanding the interplay between dynamical stability and scalar curvature rigidity of Einstein manifolds. 

Let us first recall the notions for dynamical stability and instability of Einstein metrics.
\begin{definition}
 An Einstein manifold $(M,\hat{g})$ is called dynamically stable if for any $\mathcal{U}$ of $\hat{g}$ in the space of metrics, there exists a neighbourhood $\mathcal{V}\subset\mathcal{U}$
such that for any $g_0\in\mathcal{V}$, the (appropriately normalized) Ricci flow $g_t$ starting at $g_0$  exists for all $t\geq0$ and there exists a family of diffeomorphisms $\phi_t$, $t\geq 0$ such that $\tilde{g}_t:=\phi_t^*g_t\in\mathcal{U}$ for all $t\geq0$ and
$\tilde{g}_t$
 converges to an Einstein metric in $\mathcal{U}$ as $t\to\infty$.

We call an Einstein manifold $(M,\hat{g})$ dynamically unstable if there exists a nontrivial (appropriately normalized) Ricci flow $g_t$ defined on $(-\infty,0]$ and a family of diffeomorphisms $\varphi_t$, $t\in (-\infty,0]$ such that $\phi_t^*g_t$ converges to $\hat{g}$ as $t\to-\infty$.
\end{definition}
Note that we did not fix the topology on the space of metrics in the above definition. On closed manifolds, the topology is typically $C^{k,\alpha}$ but due to parabolic estimates, one could as well just pick the $L^{\infty}$-topology. Any choice of reasonable topology will yield the same result. On the other hand, stability properties of complete noncompact manifolds may crucially depend on the decay assumptions of the metric perturbations at infinity and thus the topology has to be specified in that case.
The diffeomorphisms appearing in the above definition reflect the fact that the Ricci tensor is a diffeomorphism-invariant quantity and thus the Ricci flow is only parabolic in the strict sense after adding a suitable Lie derivative term in the evolution equation. For a particular choice of Lie derivative, this process is called de Turck trick and it yields the Ricci de Turck flow  which often converges in the strict sense.
Finally, note that one can in general not expect that the flow will converge to $\hat{g}$ as there may be  a nontrivial moduli space of Einstein metrics near $\hat{g}$.

Rigidity questions involving scalar curvature have received a new boom in the last decade, see \cite{GL22} for an overview. We will focus on Einstein metrics and use the following notion of scalar curvature rigidity in this article.
\begin{definition}\label{def:scr}
An Einstein metric $(M,\hat{g})$ is called scalar curvature rigid, if there exists a neighborhood $\mathcal{U}$ near $\hat{g}$ in the space of metrics such that there is no metric $g\in\mathcal{U}$ which agrees with $\hat{g}$ outside a compact set and satisfies $\scal_g\geq\scal_{\hat{g}}$ everywhere and $\scal_g>\scal_{\hat{g}}$ somewhere.
\end{definition}
Scalar curvature rigidity of Euclidean space is by now classical due to the rigidity part of the positive mass theorem for asymptotically Euclidean manifolds \cite{SY79,Wit81}. In fact, the rigidity is even global on the space of metrics: Any metric on $\R^n$ with $\scal\geq0$ which agrees with the Euclidean metric $g_{eucl}$ outside a compact set must be isometric to $g_{eucl}$. An analogous statement holds for hyperbolic space $\Hyp^n$ for metrics with $\scal\geq -n(n-1)$, thanks to a positive mass theorem for asymptotically hyperbolic metrics \cite{CH03,MO89}. Inspired by these examples, it was widely expected that the upper half of the round sphere is scalar curvature rigid as well but this conjecture was surprisingly disproved by Brendle-Marques-Neves \cite{BMN11}.

Concepts of mass exist also for other asymptotic geometries. For these, there exist versions of the positive mass theorem for spin manifolds which imply the scalar curvature rigidity of certain noncompact  Ricci-flat manifolds with a parallel spinor \cite{Dah97,Dai04}.

\begin{remark}
If $(M,\hat{g})$ is an Einstein manifold with nonzero scalar curvature, then it is not scalar curvature rigid, unless it admits nontrivial conformal Killing fields. In that case, we additionally impose the condition that the metric $g$ admits the same volume on the compact set $\mathrm{supp}(g-\hat{g})$.
\end{remark}

In \cite{DK24}, Dahl and the author made the attempt to systematically characterize the scalar curvature rigidity of Einstein manifolds. In the closed setting, we discovered a strong relation to dynamical stability through the interrelation of some involved geometric functionals. For open and bounded regions within an Einstein manifold, we characterize modulo some conditions scalar curvature rigidity in terms of the Dirichlet eigenvalues of the Einstein operator on the closure of that region. For complete noncompact manifolds which are asymptotically locally Euclidean (ALE), scalar curvature rigidity is strongly linked to the positive mass theorem. Motivated by this relation, Dahl, McCormick and the author 
defined a new interesting mass-like quantity for   asymptotically hyperbolic (AH) manifolds \cite{DKM23}.
Further links to dynamical stability in the ALE and the AH setting have been established over the last few years by Deruelle-Ozuch \cite{DO20,DO23} and work of the author with Deruelle \cite{DK21}, Petersen \cite{KP20} and Yudowitz \cite{KY25}, respectively.

This overview article is structured according to the different geometric settings which have been investigated. In Section \ref{sec_cpt}, we discuss the closed case and in Section \ref{sec_open}, the open case. Section \ref{sec_noncpt} discusses two cases of noncompact settings: The ALE case in Subsection \ref{subsec_ALE} and the AH case in Subsection \ref{subsec_AH}.

\vspace{3mm}
\noindent\textbf{Acknowledgements.}
We acknowledge support by the DFG priority program
SPP 2026 \emph{Geometry at Infinity}.
This survey will appear in the final volume reporting on the
SPP 2026. The author  is grateful for this opportunity.
He would further like to thank Uwe Semmelmann
for helpful and inspiring discussions.

\section{Closed manifolds}\label{sec_cpt}
A very important functional on closed manifolds is the (conformal) Yamabe invariant of a conformal class of metrics which is given by
\begin{equation*}
Y(M^n,[g]):=\inf_{\tilde{g}\in [g]} \mathrm{vol}(M,\tilde{g})^{\frac{2}{n}-1}\int_M\scal_{\tilde{g}}\dv_{\tilde{g}}.
\end{equation*}
In the following we think of the conformal Yamabe invariant as the functional $Y:\mathcal{M}\ni g\mapsto Y(M,[g])$,
where $\mathcal{M}$ is the space of smooth Riemannian metrics on the manifold $M$. For completeness, let us also recall that the smooth Yamabe invariant is defined by
\begin{align*}
Y(M):=\sup_{g\in \mathcal{M}}Y(M,[g])
\end{align*}
and that the supremum is always finite and satisfies the bound
\begin{equation*}
Y(M^n)\leq Y(S^n)=Y(S^n,[g_{rd}]).
\end{equation*}
The conformal Yamabe invariant is an important functional appering in scalar curvature rigidity problems as the theorems below show.

Recall also that a symmetric $2$-tensor is called a \emph{TT-tensor}, if it is pointwise trace-free and divergence-free. On an Einstein manifold, the space  of TT-tensors is preserved by $\Delta_E$. We call a closed Einstein manifold $(M,\hat{g})$ \emph{linearly stable}, if all eigenvalues of $\Delta_E$  on TT-tensors are nonnegative and \emph{linearly unstable} otherwise. It is called \emph{integrable}, if every  $h\in\mathrm{ker}(\Delta_E|_{TT})$ is tangent to a $C^1$-family of Einstein metrics through $\hat{g}$.

Now we are ready to formulate our main statement for closed Ricci-flat manifolds. The topology we use for the locality statements here and in the remainder of the section is the $C^{2,\alpha}$-topology on the space of metrics.
\begin{theorem}\label{thm:rigiditycptRicciflat}
Let $(M,\hat{g})$ be a closed Ricci-flat manifold. Then the following are equivalent.
\begin{enumerate}
\item[(i)] $\hat{g}$ is a local maximizer of the conformal Yamabe invariant.
\item[(ii)] Close to $\hat{g}$, there is no metric $g$ with $\scal_{g}\geq0$ and $\scal_{g}>0$ somewhere.
\item[(iii)] Close to $\hat{g}$, there is no metric of constant positive scalar curvature.
\item[(iv)] Any scalar-flat metric close to $\hat{g}$ is also Ricci-flat.
\item[(v)] $\hat{g}$ is a local maximizer of the $\lambda$-functional.
\item[(vi)] $\hat{g}$ is dynamically stable under the Ricci flow.
\end{enumerate}
Furthermore, the negation of one (and hence any) of the assertions (i)-(v) is equivalent to dynamical instability.
\end{theorem}
In this complete form, the theorem was first formulated by Dahl and the author \cite{DK24} but most implications where shown previously.
 The equivalence of the assertions (i)-(iv) follows to a large extent from the resolution of the Yamabe problem \cite{Sch84}, with additional contributions from \cite{DK24,DWW05,Koi79}.

Recall that the $\lambda$-functional associates to each metric $g$ the smallest eigenvalue of the operator $4\Delta_g+\scal_g$. It follows from the fundamental work of Perelman that $\lambda$ is a Liapunov functional for the Ricci flow. The equivalence (v)-(vi) is therefore natural to expect but still highly nontrivial to prove. It was established by Haslhofer-Mueller \cite{HM14}.

Finally, the bridge between the equivalences (i)-(iv) and (v)-(vi) is built by relating the $\lambda$-functional with the smallest eigenvalue $\lambda_{Y}$ of the Yamabe operator
\begin{align*}
4\frac{n-1}{n-2}\Delta_g+\scal_g,
\end{align*}
as observed in \cite{Kro20}.
If (i) holds, then $\hat{g}$ is a local maximizer of $\lambda_{Y}$ and the trivial inequality $\lambda\leq \lambda_{Y}$ implies (v). On the other hand, if (i) fails, we find a metric $g$ of positive scalar curvature near $\hat{g}$ which therefore satisfies $\lambda(g)>0=\lambda(\hat{g})$ and (v) fails. An argument using the second variation yields the following criterion.
\begin{theorem}\label{thm:linstabRicciflat}
Let $(M,\hat{g})$ be a closed Ricci-flat manifold. Then linear stability and integrability imply (i)-(vi) in Theorem \ref{thm:rigiditycptRicciflat}. Linear instablity implies the negation of (i)-(v) in Theorem \ref{thm:rigiditycptRicciflat} and dynamical instablity.
\end{theorem}
\begin{remark}\label{rem:cptRicciflat}
All known examples of closed Ricci-flat manifolds admit a universal cover with a parallel spinor and are therefore linearly stable and integrable, see \cite{AKWW}. Thus, all these examples satisfy assertions (i)-(vi) in Theorem \ref{thm:rigiditycptRicciflat}. 
\end{remark}

It is now natural to ask for an analogue of Theorem \ref{thm:rigiditycptRicciflat} for closed Einstein manifolds. In this case, it is natural to impose a volume constraint in order to avoid rescalings of metrics which obviously just rescale the scalar curvature.
Therefore we call a closed Einstein manifold scalar curvature rigid if Definition \ref{def:scr} holds for all $g$ near $\hat{g}$ which have the same volume.
 To formulate the theorem, let $\mathcal{M}_c$ be the set of smooth Riemannian metrics on $M$ of volume $c>0$.
\begin{theorem}[\cite{DK24,Kro20}]\label{thm:rigiditycptEinstein}
Let $(M,\hat{g})$ be a closed Einstein manifold of nonzero scalar curvature and $c=\mathrm{vol}(M,\hat{g})$. Then the following are equivalent.
\begin{enumerate}
\item[(i)] $\hat{g}$ is a local maximizer of the Yamabe invariant.
\item[(ii)] Close to
$g$, there is no metric $g\in\mathcal{M}_c$ such that $\scal_{g}\geq \scal_{\hat{g}}$ and $\scal_{g}>\scal_{\hat{g}}$ somewhere.
\item[(iii)] Close to $\hat{g}$, there is no constant scalar curvature metric $g\in\mathcal{M}_c$ such that $\scal_{g}> \scal_{\hat{g}}$.
\item[(iv)] Any metric $g\in\mathcal{M}_c$ close to $\hat{g}$ with $\scal_{g}=\scal_{\hat{g}}$ is also Einstein.
\end{enumerate}
Furthermore, the following are equivalent.
\begin{enumerate}
\item[(v)] $\hat{g}$ is a local maximizer of the shrinker entropy (if $\scal_g>0$), respectively the expander entropy (if $\scal_g<0$).
\item[(vi)] $\hat{g}$ is dynamically stable under the volume-normalized Ricci flow.
\end{enumerate}
In addition, the negation of (v) is equivalent to dynamical instablity.\\
\noindent If $\scal_g<0$, (i)-(v) are equivalent to (v)-(vi). 
If $\scal_g>0$, we have the following:
\begin{enumerate}
\item[(vii)] Assertion (i) implies (v), if we have in addition the inequality $\mathrm{spec}_+(\Delta)>2\frac{\scal}{n}$ .
\item[(viii)] Assertion (v) implies (i) and the inequality $\mathrm{spec}_+(\Delta)\geq 2\frac{\scal}{n}$. Furthermore, we have $\int_Mv^3\dv= 0$ for all $v\in C^{\infty}(M)$ with $\Delta v= 2\frac{\scal}{n}v$.
\end{enumerate}
\end{theorem}
Recall that the shrinker and expander entropies are Liapunov functionals for the respective adapted versions of the Ricci flow. Again, the equivalence (v)-(vi) is natural to expect and was established by the author in \cite{Kro20} by an argument similar to Haslhofer-M\"{u}ller in the Ricci-flat case \cite{HM14}.

The equivalence of the assertions  (i)-(iv) was shown in \cite{DK24}.
In contrast to the Ricci-flat case, none of these can be related to the resolution of the Yamabe problem. Instead, the theorem heavily relies on the local decomposition
\begin{align*}
\mathcal{M}\cong C^{\infty}_{>0}(M)\times \mathcal{C}_c,
\end{align*}
near the pair $(1,\hat{g})$, which was established by Koiso \cite{Koi79}.
Here, $\mathcal{C}_c$ is the manifold of constant scalar curvature metrics with volume $c$ and $C^{\infty}_{>0}(M)$ acts on $ \mathcal{C}_c$ by multiplying by a conformal factor. The advantage of using $\mathcal{C}_c$ is that many geometric functionals, including the shrinker and expander entropies as well as the smallest eigenvalue of any operator of the form $\alpha\Delta+\scal$ with $\alpha>0$, can explicitly be computed on metrics $g\in \mathcal{C}_c$ in terms of the (constant) scalar curvature metrics. In order to completely understand their variational behaviour, it therefore remains to compute the variations in conformal directions and the sign of these variations relies on spectral properties of the Laplace-Beltrami operator $\Delta$.
\begin{theorem}
Let $(M,\hat{g})$ be a closed Einstein manifold of nonzero scalar curvature. Then linear stability and integrability imply (i)-(iv) in Theorem \ref{thm:rigiditycptRicciflat}. If we in addition have $\mathrm{spec}_+(\Delta)>2\frac{\scal}{n}$, then (v)-(vi) hold as well.
Linear instablity implies the negation of (i)-(iv) in Theorem \ref{thm:rigiditycptRicciflat} and dynamical instablity.
\end{theorem}
\begin{remark}\label{rem:closed_negEinstein}
All known examples of closed negative Einstein manifolds satisfy the assertions (i)-(vi) in Theorem \ref{thm:rigiditycptEinstein}. 
\end{remark}
Remark \ref{rem:cptRicciflat} and Remark \ref{rem:closed_negEinstein} raise the following major open question. 
\begin{question}[Positive mass problem for closed Einstein manifolds]
Are all closed Einstein manifolds with $\scal\leq 0$ scalar curvature rigid and dynamically stable?
\end{question}
\begin{example}
For positive Einstein manifolds, there dynamically stable examples (e.g.\ $S^n$) as well as dynamically unstable ones (e.g.\ $S^n\times S^m$). It seems that dynamically stable examples are rather special and require a high degree of symmetry.

Note also that dynamical stability and scalar curvature rigidity are not equivalent for positive Einstein manifolds. Dynamical stability imples scalar curvature rigidity but the converse implcation is not true. Complex projective space $\C\mathrm{P}^{m}$ (of real dimension $n=2m$) is for $m>1$ dynamically unstable although it is scalar curvature rigid because it admits a function $v$ satisfying $\Delta v=2\frac{\scal}{n}v$ and $\int_Mv^3\dv\neq 0$. 
\end{example}

\section{Open manifolds}\label{sec_open}
Let us call an open Einstein manifold $(M,g)$ \emph{linearly stable}, if the number
\[
\mu_1(\Delta_E|_{TT},M):=\inf\left\{(\Delta_Eh,h)_{L^2}\mid h\in C^{\infty}_c(TT),\left\|h\right\|^2_{L^2}=1\right\}
\]
is nonnegative and \emph{linearly unstable} otherwise. For open Ricci-flat manifolds, we can state the following:
\begin{theorem}[\cite{DK24}]\label{thm:rigidityopenRicciflat}
Let $(M,\hat{g})$ be an open Ricci-flat manifold which does not admit a linear function, that is, there is no nonconstant function $f$ with $\nabla^2f\equiv0$. Then the following are equivalent:
\begin{enumerate}
\item[(i)] $(M,\hat{g})$ is linearly stable.
\item[(ii)] Close to $\hat{g}$, there is no metric $g$ with $g-\hat{g}|_{M\setminus K}\equiv0$ for some compact set $K\subset M$, which additionally satisfies
\[
\scal_g\geq 0, \qquad \scal_{\hat{g}}(p)>0 \text{ for some }p\in M.
\]
\item[(iii)] If $g$ is a metric close to $\hat{g}$ with $\scal_{g}\equiv0$ and $g-\hat{g}|_{M\setminus K}\equiv0$ for some compact set $K\subset M$, then $g$ is isometric to $\hat{g}$.
\end{enumerate}
\end{theorem}
\begin{remark}
The implications (i)$\Rightarrow$(ii) and (i)$\Rightarrow$(iii) do even hold in the presence of linear functions.
\end{remark}
 The typical example to which Theorem \ref{thm:rigidityopenRicciflat} applies is an open subset $\Omega$ of a complete Einstein manifold $M$ without parallel vector fields. 
By domain monotonicity of $\mu_1(\Delta_E|_{TT},\Omega)$ and the fact that $\mu_1(\Delta_E|_{TT},\Omega)\to \infty$ as $\Omega$ shrinks to a point, we see that sufficiently small open subsets of a given Ricci-flat manifold $M$ are always scalar curvature rigid.
If $M$ is linearly unstable, we will have $\mu_1(\Delta_E|_{TT},\Omega)=0$ once $\Omega$ grows to a critical size. For larger $\Omega$, we will  have $\mu_1(\Delta_E|_{TT},\Omega)<0$ and scalar curvature rigidity fails.
\begin{example}
Open Ricci-flat manifolds whose universal cover carries a parallel spinor are linearly stable \cite{DWW05} and hence scalar curvature rigid. 

In contrast to the closed case, there exist also examples of unstable open Ricci-flat manifolds, for example the Ricci-flat ALF families of the Riemannian Kerr, Taub-Bolt and Chen-Teo metrics, see \cite{BO23}.
\end{example}
\noindent The Einstein version of Theorem \ref{thm:rigidityopenRicciflat} is as follows:
\begin{theorem}[\cite{DK24}]\label{thm:rigidityopenEinstein}
Let $(M,\hat{g})$ be an open Einstein manifold satisfying the following two assumptions:
\begin{enumerate}
\item[(A)] $(M,\hat{g})$ is not locally isometric to a warped product.
\item[(B)] If $\scal_{\hat{g}}>0$, $M$ is the interior of a compact manifold $\overline{M}$ with smooth boundary whose first nonzero Neumann eigenvalue satisfies
\begin{equation} \label{eq_spectral_inequality}
\mu_1^{NM}(\Delta_{\hat{g}},\overline{M})>\frac{\scal_{\hat{g}}}{n-1}.
\end{equation}
\end{enumerate}
Then the following are equivalent:
\begin{enumerate}
\item[(i)] $(M,\hat{g})$ is linearly stable.
\item[(ii)] Close to $\hat{g}$, there is no metric $g$ with 
\[
g-\hat{g}|_{M\setminus K}\equiv0,\qquad \volume(K,g)=\volume(K,\hat{g})
\]
for some compact set $K\subset M$, which additonally satisfies
\[
\scal_g\geq\scal_{\hat{g}},\qquad \scal_g(p)>\scal_{\hat{g}}(p)
\text{ for some }p\in M.
\]
\item[(iii)] If $g$ is a metric close to $\hat{g}$ with $\scal_{g}\equiv\scal_{\hat{g}}$ and
\[
g-\hat{g}|_{M\setminus K}\equiv0,\qquad
\volume(K,g)=\volume(K,\hat{g})
\]
for some compact set $K\subset M$, then $g$ is isometric to $\hat{g}$.
\end{enumerate}
\end{theorem}
\begin{remark}
The two implications (i)$\Rightarrow$(ii) and (i)$\Rightarrow$(iii) also do hold without Assumption (A). 
Conversely, the two implications (ii)$\Rightarrow$(i) and (iii)$\Rightarrow$(i) also do hold without Assumption (B).
\end{remark}
\begin{example}
Open Einstein manifolds of special holonomy and negative scalar curvature are linearly stable and therefore satisfy (ii) and (iii) as well. Two types of special holonomy are possible for Einstein manifolds with $\scal<0$ and the linear stability of both cases has been established separately.
The K\"{a}hler-Einstein case follows from work by Koiso \cite{Koi83} (see also \cite[p.\ 361-365]{Bes08})
and the quaternion-K\"{a}hler case follows from recent work by Semmelmann and the author \cite{KS24}.

On the other hand, there arew open Einstein manifolds which are linearly unstable, for example the warped product metrics discussed in \cite{Kro17} and the AH metrics discussed later in Example \ref{ex_AH}.
\end{example}
\begin{remark}
Theorem \ref{thm:rigidityopenRicciflat} and Theorem \ref{thm:rigidityopenEinstein} can be regarded as localized versions of parts of the respective theorems in Section \ref{sec_cpt}. It is not clear yet whether characterizations in terms of dynamical stability hold as well.
\end{remark}
\section{Complete noncompact manifolds}\label{sec_noncpt}
Note that the scalar curvature rigidity results in Section \ref{sec_open} in particular apply to complete noncompact manifolds.
 Dynamical stability is harder to establish in the noncompact setting although substantial progress has been made recently. 
 We will focus on ALE manifolds, which are a typical class contaning Ricci-flat metrics and AH manifolds, which are typical for negative Einstein metrics.
 As we will see, the methods and results heavily depend on the asymptotic geometry. In both cases, we call an Einstein manifold $(M,\hat{g})$ linearly stable, if $\Delta_E|_{TT}\geq0$ in the $L^2$-sense and integrable, if any $h\in\mathrm{ker}_{L^2}(\Delta_E|_{TT})$ is tangent to a family of Einstein manifolds through $\hat{g}$.
\subsection{ALE manifolds}\label{subsec_ALE}
Recall that a complete Riemannian manifold $(M^n, g)$ is called asymptotically locally Euclidean (ALE), if there is a compact subset $K \subset M$, a real number $\tau>0$ (called the order) and a diffeomorphism $\varphi:  M \backslash K \to \R^n_{> 1}/\Gamma$ such that
\begin{equation*} \label{eq: decay metric}
	\mathring{\nabla}^k(\varphi_*g - \mathring g) \in O\left(r^{-\tau-k}\right),\qquad k\in\N_0.
\end{equation*}
Here, $\mathring{g}$ denotes the Euclidean metric.

Many examples of Ricci-flat ALE are known, for example the four-dimensional hyperkähler ALE manifolds provided by Kronheimer \cite{Kro89}. A motivation for studying dynamical stability in this context is that Ricci-flat ALE manifolds can appear as blowup limits of singularities which develop under the Ricci flow \cite{App23}. In order to know whether the singularity is generic (i.e.\ it can not be removed by a  perturbation of the initial data), it is crucial to know whether the blowup is stable as a stationary point of the Ricci flow.

Dynamical stability for $\R^n$ was first proven by Schulze,
Schn\"{u}rer and Simon \cite{SSS08}. A refinement of this result with a more precise analysis of the convergence rate was later established by Koch and Lamm \cite{KL12}. These proofs rely heavily on the explicit geometry of $\R^n$ and cannot be generalized to the ALE setting.

This geometric situation is much more challenging than the closed case, which can already be read off from spectral properties of the Lichnerowicz Laplacian $\Delta_L$. 
In the closed case, we have 
\begin{equation}\label{eq : strict positivity}
	\Delta_L|_{\mathrm{ker}^{\perp}}\geq c>0,
\end{equation}
provided that linear stability holds. On Ricci-flat ALE manifolds, the continuous spectrum of $\Delta_L$  is always $[0,\infty)$ and thus, \eqref{eq : strict positivity} can never hold in this setting.
Instead, we have the weaker inequality
\begin{equation}\label{eq : weak positivity}
	\Delta_L|_{\ker_{L^2}(\Delta_L)^{\perp}}\geq c\nabla^*\nabla>0.
\end{equation}
The first general nonlinear stability result for linearly stable and integrable Ricci-flat ALE manifolds was proven
by Deruelle and the author \cite{DK21}, where the initial data was assumed to be $L^{2}\cap L^{\infty}$-close and no convergence rate was established. Moreover, only convergence of
the Ricci-de Turck flow was shown in \cite{DK21}, but not convergence of the Ricci flow. The proof relies on a-priori estimates based on \eqref{eq : weak positivity}, which requires the $L^{2}\cap L^{\infty}$-topology.

A significant improvement of this result was established by Petersen and the author, which is as follows:
\begin{theorem}[\cite{KP20}]\label{thm:ALE_dyn_stab}
	Let $(M^n,\hat{g})$ be an ALE manifold, which carries a parallel spinor and is integrable. Then for each $q\in (1,n)$, it is dynamically stable with respect to the $L^q\cap L^{\infty}$-topology.
 Moreover, if $g_0-\hat{g}\in L^p$ for some $p\in (1,q]$ and $g_t$, $t\geq0$ is the Ricci flow starting at $g_0$, there exists a smooth family of Ricci-flat metrics $\overline{g}_t$, $t\geq0$ smoothly converging to a Ricci-flat limit $\overline{g}_{\infty}$, and a family of diffeomorphisms $\phi_t$, $t\geq0$, such that following convergence rates do hold: 
\begin{enumerate}
	\item[(i)] For each $k\in\N_0$ and $\tau>0$, there exists a constant $C=C(\tau, k, \hat g)$ such that for all $t\geq 1$, we have
	\begin{align}\label{eq: decay rate h}
		\left\|\overline{g}_t-\overline{g}_{\infty}\right\|_{C^k}\leq C\cdot t^{1-\frac{n}{p}+\tau}.
	\end{align}
	\item[(ii)] For $r\in [p,\infty]$ and $k \in\N_0$ such that $\frac{n}{2}\left(\frac{1}{p}-\frac{1}{r}\right)+\frac{k}{2}<\frac{n}{2p}$, there exists a constant $C = C(p, q, k, \hat g)$ such that for all $t\geq 1$, we have
	\begin{align}\label{eq : decay rate k 1}
	\left\|\nabla^k (\phi_t^*g_t-\overline{g}_t)\right\|_{L^{r}} \leq C\cdot t^{-\frac n2 \left(\frac1p - \frac{1}{r}\right)-\frac{k}{2}}.
	\end{align}
	\item[(iii)] For $r\in [p,\infty]$ and $k\in\N_0$ such that $\frac{n}{2}\left(\frac{1}{p}-\frac{1}{r}\right)+\frac{k}{2}\geq\frac{n}{2p}$ and for each $\tau>0$ there exists a constant $C = C(p, q,k,\tau, \hat g)$ such that for all $t\geq 1$, we have
	\begin{align}\label{eq : decay rate k 2}
	\left\|\nabla^k (\phi_t^*g_t-\overline{g}_t)\right\|_{L^{r}} \leq C\cdot t^{-\frac{n}{2p}+\tau}.
	\end{align}
\end{enumerate}
\end{theorem}
Inspired by an approach by Koch and Lamm \cite{KL12}, we established the Ricci-de Turck flow as the fixed point of a contraction map on a space of time-dependent tensors equipped with a specific norm. Thereby, we proved long-time existence and convergence at once, and the rates follow from the definition of the norm. A cruical ingredient are novel estimates for the heat kernel of $\Delta_L$ and its derivatives, which Petersen and the author developed in an earlier paper \cite{KP21}. The assumption of having a parallel spinor was important in developing these heat kernel estimates. Note that it also implies linear stablity \cite{DWW05}.

If additionally $p<\frac{3n}{4}$, we  were able to get rid of the diffeomorphisms $\phi_t$ and prove convergence of the Ricci flow directly. Furthermore, we could use Theorem \ref{thm:ALE_dyn_stab} the following scalar curvature rigidity result:
\begin{theorem}[\cite{KP20}]\label{thm:ALE_scr}
	Let $(M^n,\hat{h})$ be an ALE Ricci-flat spin manifold which is integrable and carries a parallel spinor. Then for each $q\in (1,n)$, there exists a $L^q\cap L^{\infty}$-neighbourhood $\mathcal{U}$ of $\hat{h}$ in the space of metrics such that each smooth metric $g\in\mathcal{U}$  on $M$ satisfying
	\begin{equation*}
	 \scal_g\geq0,
	 	 \qquad\text{ and }\qquad\left\|g-\hat{h}\right\|_{L^p}<\infty
	\end{equation*}
	for some $p\in \left[1,\frac{n}{n-2}\right)$ is Ricci-flat.
\end{theorem}
The proof is based on the following contradiction argument: Suppose we have a metric $g$ as in the theorem. Then by a heat equation, which the scalar curvature satsifies along the Ricci flow, $\left\|\scal_{g_t}\right\|_{C^0}\geq C t^{-\frac{n}{2}}$ unless it is Ricci-flat. On the other hand, using the proof of Theorem \ref{thm:ALE_dyn_stab} , one can show that $\left\|\scal_{g_t}\right\|_{C^0}\leq C t^{-\frac{n}{2p}-1}$. This leads to a contradiction for $p<\frac{n}{n-2}$ unless $g$ is Ricci-flat.

A different approach towards the dynamical stability of Ricci-flat ALE manifolds has been developed by Deruelle and Ozuch \cite{DO20,DO23}, who use an ALE version of the $\lambda$-functional. The first attempt to define such a functional was made by Haslhofer \cite{Has11}, who set it to be
\begin{align*}
\lambda^0_{ALE}(g)=\inf_{\omega-1\in C^{\infty}_{c}}\int_M(4|\nabla\omega|_g+\mathrm{scal}_g\omega^2)\,\mathrm{dV}_g.
\end{align*}
Deruelle and Ozuch realized that it is very convenient to incooperate the ADM-mass
\begin{align*}
m_{ADM}(g):=\lim_{R\to\infty}m_{ADM}(g,R):=\lim_{r\to\infty}\int_{\partial B_R}\langle \mathrm{div}_{\mathring{g}}\varphi_*g-d\mathrm{tr}_{\mathring{g}}\varphi_*g,\nu\rangle\,\mathrm{dS}_{\mathring{g}}
\end{align*}
into the functional and defined it as
\begin{align*}
\lambda_{ALE}(g)=\inf_{\omega-1\in C^{\infty}_{c}}\lim_{r\to\infty}
\left(\int_{B_R}
(4|\nabla\omega|_g+\mathrm{scal}_g\omega^2)
\,\mathrm{dV}_g-m_{ADM}(g,r)\right).
\end{align*}
This refinement has the advantage that it is well-defined for a larger class of metrics and, in contrast to $\lambda^0_{ALE}$, does not require $\scal_g$ to be integrable. For $\tau\in (2-n,\frac{2-n}{2})$, it is an analytic functional on a $C^{2,\alpha}_{\tau}$-neighborhood of a Ricci-flat ALE metric. The corresponding weighted H\"{o}lder norm is defined by
\begin{equation*}
\left\|u\right\|_{C^{2,\alpha}_{\tau}}=
\sup_M r^{-\tau}\left(\sum_{k=0}^2 r^{k}|\nabla^ku|+r^{2+\alpha}|\nabla^2u|_{C^{0,\alpha}}\right),
\end{equation*}
where $r$ is any smooth function on $M$ which on $M\setminus K$ agrees with pullback of the radius function on $\R^n$ by the diffeomorphism $\varphi$.
Under the assumption of integrability, Deruelle and Ozuch established in \cite{DO20} a Lojasziewicz-Simon inequality for $\lambda_{ALE}$ for dimensions $n\geq5$. They used it to prove dynamical stability and instability results for integrable Ricci-flat ALE manifolds with respect to the $C^{2,\alpha}_{\tau}$-topology \cite{DO23}. 

Comparing their result to Theorem \ref{thm:ALE_dyn_stab}, Deruelle and Ozuch do not assume a parallel spinor and they are able to prove an instability result as well. On the other hand, Theorem \ref{thm:ALE_dyn_stab} also works for $n=4$, it requires much weaker falloff conditions on the initial data and the convergence rates are explicitly given. 

\begin{remark}
All known examples of Ricci-flat ALE metrics carry a parallel spinor. ALE K\"{a}hler metrics are additionally known to be integrable \cite{DK21,Lund19}. Hence, these (which in particular include the Kronheimer examples) satisfy the conditions of Theorem \ref{thm:ALE_dyn_stab} and Theorem \ref{thm:ALE_scr}.
\end{remark}
Interestingly, $\lambda_{ALE}$ creates a link to the positive mass theorem (PMT). If the Ricci-flat ALE metric $\hat{g}$ is a local maximum of $\lambda_{ALE}$, then for any metric $g$ near $\hat{g}$ with $\scal_g\geq0$ and $\scal\in L^1$, we have that
\begin{equation*}
0=\lambda_{ALE}(\hat{g})\geq \lambda_{ALE}(g)=\lambda^0_{ALE}(g)-m_{ADM}(g),
\end{equation*}
and therefore, since $\scal_g\geq0$, 
\begin{equation*}
m_{ADM}(g)\geq \lambda^0_{ALE}(g)\geq 0.
\end{equation*}
This means that the assertion of the PMT  holds for metrics near $\hat{g}$. It is not hard to see that a PMT for metrics near $\hat{g}$ also implies that $\hat{g}$ is a local maximum of $\lambda_{ALE}$, although it is not directly stated in \cite{DO20}.

 A connection between an appropriate version of the $\lambda$-functional for ALE manifolds and the positive mass theorem had been already conjectured by Ilmanen before and was partly proven by Hall-Haslhofer-Siepman \cite{HHS14}. All these points discussed here raise the following question, whose resolution would form an ALE analogue of Theorem \ref{thm:ALE_dyn_stab}.
\begin{question}\label{question:Ilmanen_conj}
Is the dynamical stability of a Ricci-flat manifold $(M,\hat{g})$ equivalent to a PMT for metrics near $\hat{g}$?
\end{question}

\subsection{AH manifolds}\label{subsec_AH}
Recall that a complete Riemannian manifold $(M,g)$ is called \emph{conformally compact} of class $C^{k,\alpha}$ if $M$ is the interior of a manifold $\overline{M}$ with boundary $\partial M$ and the metric $h=\rho^2g$ can be extended to a $C^{k,\alpha}$-metric on all of $\overline{M}$. Here, $\rho:\overline{M}\to [0,\infty)$ is a boundary defining function for $\overline{M}$, that is $\rho^{-1}(0)=\partial M$ and $d\rho$ is nowhere vanishing on $\partial M$. The smooth function $r=-\log(\rho):M\to \R$ is near $\partial M$ equivalent to the distance function from a fixed point.

Note that the metric $h|_{\partial M}$ on $\partial M$ depends on $\rho$ but its conformal class is independent of it. We call $[h|_{\partial N}]$ the \emph{conformal boundary} of $(M,g)$. Because the sectional curvature $K_g$ converges to $-|d\rho|_h$ at $\partial M$, we call $(M,g)$ \emph{asymptotically hyperbolic} (AH), if $|d\rho|_h\equiv 1$ on $\partial M$.

Note further that a conformally compact Einstein metric $(M^n,g)$ has  negative Ricci curvature. Normalizing it so that $\ric_g=-(n-1)g$, we see that $(M^n,g)$ is necessarily AH. We call such manifolds \emph{Poincar\'{e}-Einstein} (PE). 

The Poincar\'{e} ball model shows that $(\Hyp^n,g_{hyp})$ is AH (and hence PE) and its conformal boundary is the round sphere. A classical result of Graham and Lee \cite{GL91} demonstrates that any conformal class near the round sphere is the conformal boundary of a PE metric on $\R^n$. This provides many examples of PE metrics and their dynamical stability under the Ricci flow is of natural interest. On the other hand, PE metrics also play an important role in the AdS-CFT correspondence and the Lorentzian cone of a PE metric is a globally hyperbolic solution of Einstein's vacuum equation. This gives plenty of motivation to study scalar curvature rigidity and positive mass problems in the context of PE manifolds and Question \ref{question:Ilmanen_conj} seeks for a connection to dynamical stability.

As a first example, the dynamical stability of $(\Hyp^n,g_{hyp})$ was established by Schulze-Schn\"{u}rer-Simon \cite{SSS11}, a result which was refined and generalized to symmetric spaces of noncompact type by Bamler \cite{Bam15}. In \cite{HJS16}, Hu-Ji-Shi proved the dynamical stability of a strictly linearly stable PE manifold. 

Inspired by recent progress in the ALE case, Dahl, McCormick and the author defined a new mass-like quantity for AH manifolds and an AH version of the expander entropy in a recent preprint \cite{DKM23}. From now on, we assume that all AH are conformally compact of class at least $C^{k,\alpha}$ with $k\geq2$.
We define the volume-renormalized mass of the AH metric $g$ (with respect to the reference AH metric $\hat{g}$) as
\begin{equation*}
m_{VR,\hat{g}}(g):=
\lim_{R\to\infty}\left(m_{ADM,\hat{g}}(g,r)+2(n-1)\mathrm{RV}_{\hat{g}}(g)\right),
\end{equation*}
where
\begin{align*}
m_{ADM,\hat{g}}(g,r)&:=\int_{r=R}\langle \mathrm{div}_{\hat{g}}\varphi_*g-d\mathrm{tr}_{\hat{g}}\varphi_*g,\nu\rangle\,\mathrm{dS}_{\hat{g}},
\\
\mathrm{RV}_{\hat{g}}(g)&:=\int_{r\leq R}\dv_g-\int_{r\leq R}\dv_{\hat{g}}.
\end{align*}
As for the ADM-mass (where we may think of the Euclidean metric $\mathring{g}$ as the reference metric), $g$ and $\hat{g}$ are in general not defined on the same manifold, but the respective AH ends need to be diffeomorphic. Moreover, the radius functions on both manifolds need to agree via the diffeomorphism $\varphi$. In this context, we say that an AH manifold $(M,g)$ is asymptotic to $(\hat{M},\hat{g})$ of order $\tau>0$, if there exist compact subsets $K\subset M$ and $\hat{K}\subset\hat{M}$ and a diffeomorphism $\varphi:M\setminus K\to \hat{M}\setminus \hat{K}$ such that
\begin{equation*}
\varphi_*g-\hat{g}\in C^{k,\alpha}_{\tau},
\end{equation*}
where the weighted H\"{o}lder space is induced by the norm
\begin{equation*}
\left\|u\right\|_{C^{k,\alpha}_{\tau}}=\left\|e^{\tau r}u\right\|_{C^{k,\alpha}}.
\end{equation*}
\begin{theorem}[\cite{DKM23}]\label{thm:welldef_mass}
Let $\tau>\frac{n-1}{2}$ and $(\hat{M},\hat{g})$ be an AH manifold which satisfies
\begin{align*}
\ric_{\hat{g}}+(n-1)\hat{g}\in C^{0,\alpha}_{\tau},
\qquad \scal_{\hat{g}}+n(n-1)\in L^1.
\end{align*}
Then, if $(M,g)$ is an AH manifold with $\scal_{{g}}+n(n-1)\in L^1$ which is asymptotic to $(\hat{M},\hat{g})$ of order $\tau$,  $m_{VR,\hat{g}}$ is finite and under a generic condition on the conformal boundary independent of the choice of $\varphi$.
\end{theorem}
Note that if $\tau>n-1$, the boundary term vanishes in the limit and the volume-renormalized mass reduces to a multiple of the renormalized volume
\begin{align*}
RV_{\hat{g}}(g)=\lim_{r\to\infty}\left(\int_{r\leq R}\dv_g-\int_{r\leq R}\dv_{\hat{g}}\right).
\end{align*}
\begin{remark}We have the additivity propery $m_{VR,\hat{g}}(g)=m_{VR,\overline{g}}(g)+m_{VR,\hat{g}}(\overline{g})$.
Therefore, the functional $g\mapsto m_{VR,\hat{g}}(g)$ shifts by a constant of one changes the reference metric and thus, not every choice of $\hat{g}$ can lead to a positive mass theorem. A PE metric seems to be a natural choice for $\hat{g}$. 
\end{remark}
There have been other mass invariants previously defined for asymptotically hyperboloidal manifolds \cite{CN01,Wan01}. The
quantity $m_{VR,\hat{g}}$ seems to be significant for two reasons: Using it as a normalization, we could in \cite{DKM23} for the first time establish the Einstein-Hilbert action for asymptotically hyperbolic manifolds in a mathematically clean way. 
Moreover, we showed that $m_{VR}$ can also be motivated from the Hamiltonian perspective on general relativity applied to expanding spacetimes asymptotic to the Milne model  $(\R_+ \times \R^n , -dt^2 + t^2 g_{hyp} )$ (see \cite{DKM25}), analogous to the Hamiltionian viewpoint for the ADM-mass for spacetimes which are asymptotically Minkowski. 

In \cite{DKM23}, we showed that $m_{VR}$ shares some analogies with $m_{ADM}$. For two-dimensional manifolds asymptotic to the hyperbolic plane (with a possible conical singularity at the origin), there is a formula for $m_{VR}$ in terms of the Euler characteristic and the cone angle which is completely analogous to a formula for $m_{ADM}$ for asymptotically conical surfaces. In all dimensions, we show the following: If $g$ is an AH metric with $\scal_g\geq -n(n-1)$, $v\in C^{2,\alpha}_{\tau}$ and $\overline{g}=e^{2v}g$ satisfies $\scal_{\overline{g}}=-n(n-1)$, we have the inequality $m_{VR,\hat{g}}(g)\geq  m_{VR,\hat{g}}(\overline{g})$, for which there is also an analogue involving $m_{ADM}$.

A positive mass theorem for $m_{VR}$ has been shown in \cite{DKM23} for AH metrics on $\Hyp^3$ asymptotic to the hyperbolic metric. In an upcoming preprint, Oronzio, Pinoy and the author are extending this result as follows:
\begin{theorem}[\cite{KOP25}]\label{thm:PMT}
Let $(M,g)$ be a three-dimensional orientable AH manifold which is asymptotic to $(\Hyp^3,g_{hyp})$ of order $\tau>1$. Assume that  $\scal_g\geq -6$ and that the second homology of $M$ does not contain any spherical classes. Then,
\begin{align*}
m_{VR,g_{hyp}}(g)\geq0
\end{align*}
and equality holds if and only if $(M,g)$ is isometric to $(\Hyp^3,g_{hyp})$.
\end{theorem}
To connect $m_{VR,\hat{g}}$ to the Ricci flow, a version of the expander entropy for AH manifolds was defined in \cite{DKM23} as
\begin{equation*}
\mu_{AH,\hat{g}}(g)=\inf_{\omega-1\in C^{\infty}_c}\mathcal{W}_{{AH},\hat{g}}(g,\omega),
\end{equation*}
where
\begin{align*}
{\mathcal{W}}_{{AH},\hat{g}}(g,\omega) 
&= 
\lim_{R \to \infty}
\Big( \int_{B_R}\Big( 4|\nabla \omega|^2 +(\scal_g+n(n-1))\omega^2\\
&\qquad\qquad +2(n-1)[(\log(\omega^2)-1)\omega^2+1]
\Big) dv_g\\
&\qquad\qquad -\left( m_{{\rm ADM},\hat{g}}(g,R)+2(n-1)RV_{\hat{g}}(g,R) \right)\Big) .
\end{align*}
Here, $\hat{g}$ is an AH reference metric satisfying the same conditions as in Theorem \ref{thm:welldef_mass}. If $g$ is asymptotic to $\hat{g}$ of order $\tau\in (\frac{n-1}{2},n-1)$, there is a unique minimizer $\omega_g$ realizing $\mu_{AH,\hat{g}}(g)$ and it depends analytically on the metric with respect to the $C^{2,\alpha}_{\tau}$-topology. The functional $g\mapsto \mu_{AH,\hat{g}}(g)$ is monotonically increasing under the AH Ricci flow $\partial_tg=-2\ric_{g}-2(n-1)g$ and constant only along a constant flow of a PE metric.
Theorem \ref{thm:AH_Ilmanen_conjecture} below provides a satisfying analogue of Theorem \ref{thm:rigiditycptRicciflat} and Theorem \ref{thm:rigiditycptEinstein} for PE manifolds which can also be regarded as a resolution of an AH version of Ilmanens conjecture for ALE manifolds. The topology used here is the $C^{2,\alpha}_{\tau}$ topology for $\tau\in (\frac{n-1}{2},n-1)$.
\begin{theorem}[\cite{DKM23,KY25}]\label{thm:AH_Ilmanen_conjecture}
Let $(M,\hat{g})$ be a PE manifold. Then the following are equivalent:
\begin{enumerate}
\item[(i)] $m_{VR,\hat{g}}(g)\geq 0$ for all metrics $g$ near $\hat{g}$ with $\scal_{g}\geq -n(n-1)$.
\item[(ii)] $RV_{\hat{g}}(g)\geq 0$ for all metrics $g$ near $\hat{g}$ with $\scal_{g}\geq -n(n-1)$ which additionally satisfy $g-\hat{g}=O(e^{-(n-1+\epsilon)r})$ for some $\epsilon>0$.
\item[(iii)] $\hat{g}$ is a local maximizer of $\mu_{AH,\hat{g}}$.
\item[(iv)]  $\hat{g}$ is $L^2\cap L^{\infty}$-dynamically stable.
\end{enumerate}
The negations of the assertions (i)-(iii) are equivalent to dynamical instability.
\end{theorem}
The implication (iv)$\Rightarrow$(ii) was already established by Hu-Ji-Shi in \cite{HJS16}. It is very interesting to see that property (ii) is strong enough to imply dynamical stability.

Let us discuss some aspects of the proof. The equivalence of (i) and (ii) essentially follows from a density argument. Inspired by the closed case, the equivalence of (i) and (ii) uses a local decomposition of the space of metrics which is as follows: There is a neighborhood $\mathcal{U}$ of $\hat{g}$ in the space of metrics and a neighborhood $\mathcal{V}$ of $1$ in the space of functions such that multiplication leads to a diffeomorphism
\begin{equation}
\mathcal{U}\cong \mathcal{V}\times \left\{g\in \mathcal{U}\mid \scal_g=-n(n-1)\right\}. 
\end{equation}
Furthermore, we use the observation that if $\scal_g=-n(n-1)$, we have $\omega_g\equiv 1$ and $\mu_{AH,\hat{g}}(g)=-m_{VR,\hat{g}}(g)$. Therefore, it remains to consider the variations of mass and entropy in conformal directions. As in the closed case, the equivalence of (iii) and (iv) follows from a suitable Lojasciewicz-Simon inequality.
Surprisingly, the proof of this inequality as well as the proofs of dynamical stability and instability are up to minor modifications the same as in the closed case. In particular, establishing dynamical stability is substantially simpler than in the ALE case. This seems to rely on two crucial properties of the operator $\Delta_E$ which are as in the compact case. As an operator $\Delta_E:H^k\to H^{k-2}$ between unweighted Sobolev spaces, it is Fredholm and
$\Delta_E|_{\mathrm{ker}}\geq c>0$ in the case of linear stability, because
 the continuous part of its spectrum is strictly positive.

\begin{example}\label{ex_AH}
Hyperbolic space is known to be dynamically stable. The familes of the Riemannian Anti-de Sitter versions of the Schwarzschild, Taub-Nut and Taub-Bolt metrics provide a large class of examples, of which some are dynamically stable and some are dynamically unstable \cite[Example 1.7]{KY25}. There als some threshold values for the involved parameters, for which it is not clear whether the corresponding metric is dynamically stable or unstable.
\end{example}



\end{sloppypar}
\end{document}